# NONSTANDARD LIMIT THEOREM FOR INFINITE VARIANCE FUNCTIONALS

By Allan Sly and Chris Heyde

*University of California, Berkeley, Australian National University and Columbia University*

We consider functionals of long-range dependent Gaussian sequences with infinite variance and obtain nonstandard limit theorems. When the long-range dependence is strong enough, the limit is a Hermite process, while for weaker long-range dependence, the limit is $\alpha$-stable Lévy motion. For the critical value of the long-range dependence parameter, the limit is a sum of a Hermite process and $\alpha$-stable Lévy motion.

**1. Introduction.** Define the stationary Gaussian sequence

$$X_i = \sum_{j=0}^{\infty} b_j \xi_{i-j},$$

where the $\xi_i$ are independent standard Gaussians and $b_j = j^{H-3/2} L_1(j)$, with $\frac{1}{2} < H < 1$, $L_1(i)$ slowly varying and $\sum_{j=0}^{\infty} b_j^2 = 1$. The $\{X_i\}$ are then long-range dependent with Hurst index $H$. Denote the filtration generated by $(\xi_i)_{i \in \mathbb{Z}}$ as $(\mathcal{F}_i)$. The Hermite polynomials form an orthogonal basis for $L^2(\mathbb{R}, e^{-x^2/2})$, so when $f$ is a function such that $Ef(X_1)^2 < \infty$, the chaos decomposition of $f(X_i)$ is given by

$$(1.1) \qquad f(X_i) = \sum_{k=0}^{\infty} f_k h_k(X_i),$$

where the sum is interpreted by convergence in $L^2$ and $f_k = \frac{1}{k!} E(f(X_1) h_k(X_1))$. The Hermite rank of $f$ is the smallest $k \geq 1$ such that $f_k \neq 0$. For more details on Hermite polynomials and their relationship to Gaussian Hilbert spaces, see [8]. The limits of normalized sums of the form $\sum_{i=1}^{\lfloor nt \rfloor} f(X_i) - Ef(x_i)$ were established in [5, 14] and [2], when $Ef(X_1)^2 < \infty$, in terms of the Hermite rank $\kappa$ and $H$. When $1 - \kappa(1-H) > 1/2$, the limit is a Hermite process, but when $1 - \kappa(1-H) \leq 1/2$, the limit is standard Brownian motion.

This paper considers the limits of normalized sums when the $f(X_i)$ have power-tailed distributions with index $0 < \alpha < 2$. These functions $f$ satisfy

$$P(f(X) > x) \sim \frac{1+\beta}{2} L_2(x) x^{-\alpha},$$

(1.2)

$$P(f(X) < -x) \sim \frac{1-\beta}{2} L_2(x) x^{-\alpha}$$

as $x \to \infty$, $L_2$ slowly varying and with $\beta \in [-1,1]$. We can find constants $a_n = n^{1/\alpha} L_3(n)$ such that $P(|f(X)| > a_n) \sim n^{-1}$ and $L_3$ is slowly varying. We will focus on the case when $1 < \alpha < 2$ and $f$ is centered so that $E(f(X)) = 0$.

For background, a random variable $X$ is *stable* if $cX \stackrel{d}{=} X_1 + X_2$ for some constant $c$ and where $X_1$ and $X_2$ are independent copies of $X$. Such distributions are completely categorized and the non-Gaussian stable laws have characteristic function

$$E\exp(itX) = \begin{cases} \exp\left[-\sigma^\alpha |\theta|^\alpha \left(1 - i\beta(\text{sign}(\theta))\tan\frac{\pi\alpha}{2}\right) + i\mu\theta\right], & \alpha \neq 1, \\ \exp[-\sigma^\alpha |\theta|^\alpha (1 + i\beta(\text{sign}(\theta))\ln|\theta|) + i\mu\theta], & \alpha = 1, \end{cases}$$

for $0 < \alpha < 2$, $\beta \in [-1,1]$, $\sigma > 0$ and $\mu \in \mathbb{R}$. We will denote this distribution by $S_\alpha(\sigma, \beta, \mu)$, following [11]. A process is an $\alpha$-*stable Lévy motion* if it is a Lévy process with increments which are stable; see [11] for more details.

The standard approach which relies on the chaos decomposition does not apply here. We use the hypercontractivity of Gaussian Hilbert spaces to establish that the extreme values of such processes are asymptotically independent. The limit can be either a Hermite process, $\alpha$-stable Lévy motion or a sum of both. This striking result, where normed sums of power tailed variables with infinite variance may converge to a self-similar process with finite moments of all orders, highlights the complexity associated with the domains of attraction of the Hermite processes, all of which have this property. Little is known about the domains of attraction beyond the case of fractional Brownian motion ($\kappa = 1$). Note that, in the case of i.i.d. variables (the classical central limit theorem), the only self-similar processes which appear as limits are Brownian motion and the $\alpha$-stable Lévy motions for $0 < \alpha < 2$.

Both infinite variance power tails and long-range dependence lead to partial sums increasing faster than $O(\sqrt{n})$. The different limit processes reflect the relative importance of these two effects. The composite process occurs only when the effects are balanced.

When $E|f(X)|^p < \infty$ for some $p > 1$, we can still define the $f_k$ and equation (1.1) can be interpreted as a stochastic distribution; see [7], Corollary 2.3.8. We give the Hermite rank its natural extension to $L^p$.

THEOREM 1.1. *If $1-\kappa(1-H) > \frac{1}{\alpha}$ or if $1-\kappa(1-H) = \frac{1}{\alpha}$ and $\lim_n \frac{L_1(n)^\kappa}{L_3(n)} = \infty$, then*

$$L(n)^{-\kappa} n^{-(1-\kappa(1-H))} \sum_{i=1}^{\lfloor nt \rfloor} f(X_i) \to f_\kappa R_{\kappa,H}(t),$$

*where $R_{\kappa,H}$ is the $\kappa$th Hermite process given by the multiple stochastic integral*

$$R_{\kappa,H}(t) = \int_{\mathbb{R}^\kappa} \int_0^t \prod_{i=1}^\kappa I(s > x_i) |s - x_i|^{H-3/2} \, ds \, dB(x_1) \cdots dB(x_\kappa),$$

*$B(t)$ being standard Brownian motion. If $1 - \kappa(1-H) < \frac{1}{\alpha}$, or if $1 - \kappa(1-H) = \frac{1}{\alpha}$ and $\lim_n \frac{L_1(n)^\kappa}{L_3(n)} = 0$, then*

$$a_n^{-1} \sum_{i=1}^{\lfloor nt \rfloor} f(X_i) \to R^*(t),$$

*where $R^*(t)$ is $\alpha$-stable Lévy motion with $R^*(1) \stackrel{d}{=} S_\alpha((\frac{\Gamma(2-\alpha)\cos(\pi\alpha/2)}{1-\alpha})^{1/\alpha}, \beta, 0)$. Finally, if $1 - \kappa(1-H) = \frac{1}{\alpha}$ and $\lim_n \frac{L_1(n)^\kappa}{L_3(n)} = \lambda \in (0, \infty)$, then*

$$a_n^{-1} \sum_{i=1}^{\lfloor nt \rfloor} f(X_i) \to \lambda f_\kappa R_{\kappa,H}(t) + R^*(t),$$

*where $R_{\kappa,H}$ and $R^*(t)$ are independent. Convergence is taken to mean weak convergence on $D[0,T]$ in the Skorohod $J_1$-topology.*

One might have expected the process $f(X_i)$ to be in the domain of attraction of a fractional stable process. However, it does not exhibit clustering of extreme values, as can be seen by a calculation of its extremal index. Our results imply that

$$\lim_{n \to \infty} P\left[\max_{1 \le i \le n} f(X_i) > a_{\lfloor n\tau \rfloor}\right] = e^{-\tau},$$

which shows that the extremal index is 1. On the other hand, a continuous version of Proposition 2.1 from [3] implies that fractional stable motion has extremal index strictly less than 1.

The situation when $0 < \alpha \le 1$ is simpler as only convergence to $\alpha$-stable Lévy motion is possible. While not explicitly mentioned, the one-dimensional version of the case $0 < \alpha < 1$ can be shown to follow from Lemma 5 and 6 of [4].

THEOREM 1.2.   *If $0 < \alpha < 1$, then*

$$a_n^{-1} \sum_{i=1}^{\lfloor nt \rfloor} f(X_i) \to R^*(t),$$

*while if $\alpha = 1$, then*

$$a_n^{-1} \sum_{i=1}^{\lfloor nt \rfloor} [f(X_i) - E(f(X_i)I(|f(X_i)| \leq a_n))] - \frac{2\psi t}{\pi} \to R^*(t),$$

*where $\psi = \ln \pi + \int_0^\infty u^{-2}(\sin u - u1(u \leq 1)) \, du$. Again, $R^*(t)$ is $\alpha$-stable Lévy motion with $R^*(1) \stackrel{d}{=} S_\alpha((\frac{\Gamma(2-\alpha)\cos(\pi\alpha/2)}{1-\alpha})^{1/\alpha}, \beta, 0)$ when $0 < \alpha < 1$ and $R^*(1) \stackrel{d}{=} S_1(\pi/2, \beta, 0)$ when $\alpha = 1$ and, again, convergence is in finite-dimensional distributions.*

The situation when $0 < H \leq \frac{1}{2}$ is simpler with convergence to the stable limit in Theorems 1.1 and 1.2.

The results herein remain true when we take any Gaussian sequence in the domain of attraction of fractional Brownian motion. However, the proofs involve much more tedious technical details and are no more informative, so we have confined our attention to the slightly less general result.

There are other results in the literature where functionals of long-range dependent processes have been shown to have both Gaussian and $\alpha$-stable Lévy limits depending on the parameter of long-range dependence. For instance, Surgailis (in [13]) found this behavior for certain bounded functionals of long-range dependent moving averages of heavy-tailed random variables and (in [12]) for the empirical process of another moving average of heavy-tailed random variables.

As a simple example, consider the case $f(x) = |x|^r - E(|X_i|^r)$. Proposition 3 of [6] showed that if $\frac{3}{4} \leq H < 1$ and $r > 0$, then

$$(1.3) \qquad n^{1-2H} \sum_{i=1}^{\lfloor nt \rfloor} f(X_i) \stackrel{d}{\to} f_2 R_{2,2H-1}(1)$$

as $n \to \infty$. With our results, we can extend this to all $r$. By Theorem 1.1, equation (1.3) holds when $1/(1-2H) < r < 0$, but

$$\left(\frac{2}{\pi}\right)^{r/2} n^{1-2H} \sum_{i=1}^{\lfloor nt \rfloor} f(X_i) \stackrel{d}{\to} R^*(t)$$

when $-1 < r < 1/(1-2H)$ and

$$n^r \sum_{i=1}^{\lfloor nt \rfloor} f(X_i) \stackrel{d}{\to} \left(\frac{\pi}{2}\right)^{r/2} R^*(t) + f_2 R_{2,2H-1}(1)$$

when $r = 1/(1 - 2H)$. Also, using Theorem 1.2, when $r = -1$,

$$\left(\frac{2}{\pi}\right)^{r/2} n^r \sum_{i=1}^{\lfloor nt \rfloor} \left[|X_i|^r - E\left(|X_i|^r I\left(|X_i|^r \leq \left(\frac{\pi}{2}\right)^{r/2} n^{-r}\right)\right)\right] - \frac{2\psi t}{\pi} \to R^*(t)$$

and when $r < -1$,

$$\left(\frac{2}{\pi}\right)^{r/2} n^r \sum_{i=1}^{\lfloor nt \rfloor} |X_i|^r \xrightarrow{d} CR^*(t).$$

**2. Proof.** Fix a $d$ such that $\sum_{j=d}^{\infty} b_j^2 = \theta < \min\{\alpha, 1/H\} - 1$. Then, by Theorem 5.1 of [8], the map $f(X_i) \mapsto E(f(X_i)|\mathcal{F}_{i-d})$ is a hypercontraction from $L^{1+\theta}$ to $L^2$, so $E(E(f(X_i)|\mathcal{F}_{i-d}))^2 < \infty$. With $U_i = \theta^{-1/2} \sum_{j=d}^{\infty} b_j \xi_{i-j}$, and expanding $h_k(X_i)$ according to Lemma D.1 of [7], we have $E(f(X_i)|\mathcal{F}_{i-d}) = \sum_{k=0}^{\infty} f_k \theta^{k/2} h_k(U_i)$. Then, by Theorem 5.6 of [14], when $1 - \kappa(1 - H) > 1/2$,

$$(2.1) \qquad L(n)^{-\kappa} n^{-(1-\kappa(1-H))} \sum_{i=1}^{\lfloor nt \rfloor} E(f(X_i)|\mathcal{F}_{i-d}) \to f_\kappa R_{\kappa,H}(t)$$

weakly as $n \to \infty$.

For our asymptotic independence result, we map our Gaussian sequence according to

$$X_i = X_i^{(n)} = \sum_{j=0}^{\infty} n^{1/2} b_j \left(B\left(\frac{i-j+1}{n}\right) - B\left(\frac{i-j}{n}\right)\right),$$

where $B(t)$ is some fixed Brownian motion. Brownian scaling guarantees that the distributions of the $X_i$'s do not depend on $n$. With this definition, it follows from Lemma 4.5 of [14] that the convergence in equation (2.1) can be taken as convergence in $L^2$ pointwise in $t$. We will routinely suppress the dependence on $n$ of various objects.

Fix a $c > 0$ and let $\nu_n$ be the simple point process on $(\mathbb{R}\setminus(-c,c)) \times \mathbb{R}$ given by point masses at points $(f(X_i)/a_n, t/n)$. Let $\nu$ be a Poisson point process with parameter measure

$$\nu'(dx, dy) = \begin{cases} \alpha \dfrac{1+\beta}{2} |y|^{-\alpha-1} \, dx \, dy, & y > 0, \\ \alpha \dfrac{1-\beta}{2} |y|^{-\alpha-1} \, dx \, dy, & y < 0. \end{cases}$$

LEMMA 2.1. *Let $\Delta$ be a finite union of finite intervals in $\mathbb{R}\setminus(-c,c)$ and let $Z_i = I[a_n^{-1} f(X_i) \in \Delta]$. Let $V$ be a finite variance random variable measurable in the $\sigma$-algebra generated by $B(t)$. Then,*

$$(2.2) \qquad E\left|E\left[V \sum_{i=1}^{\lfloor nt \rfloor} (Z_i - EZ_i) \,\bigg|\, \mathcal{F}_0\right]\right| \to 0$$

*as $n \to \infty$ and*

$$E \sum_{i=1}^{\lfloor nt \rfloor} \sum_{j=i+1}^{\lfloor nt \rfloor} Z_i Z_j \to \tfrac{1}{2} \mu_\Delta^2 t^2 \tag{2.3}$$

*as $n \to \infty$, where $\mu_\Delta = \nu'(\Delta \times [0,1))$.*

PROOF. It is sufficient to prove the result when $V \in H^{:l:}$, the homogeneous chaos of order $l$ of the Gaussian Hilbert space generated by $B(t)$ (see [8]). By equation (1.2),

$$nEZ_i \to \mu_\Delta \tag{2.4}$$

as $n \to \infty$. $Z_i$ can be written as $\sum_{k=0}^\infty g_k h_k(X_i)$. Applying [8], Theorem 5.1 with the map $X_i \mapsto (\ln n)^{-1} X_i$, we then obtain $\sum_{i=0}^\infty (\ln n)^{-2k} g_k^2 k! \leq \|Z_i\|_{1+(\ln n)^{-2}}^2 = (EZ_i)^{2/(1+(\ln n)^{-2})} = O(n^{-2})$ and so $g_k^2 k! = O((\ln n)^{2k} n^{-2})$. Hence, since $E(\sum_{i=1}^{\lfloor nt \rfloor} h_k(X_i))^2$ is $O(L(n)^{2\kappa} n^{2(1-\kappa(1-H))})$ when $1 - \kappa(1 - H) \geq 1/2$ and $O(n)$ when $1 - \kappa(1-H) < 1/2$, we have

$$E\left(\sum_{i=1}^{\lfloor nt \rfloor} \sum_{k=1}^{\rho} g_k h_k(X_i)\right)^2 = \sum_{k=1}^{\rho} g_k E\left(\sum_{i=1}^{\lfloor nt \rfloor} h_k(X_i)\right)^2$$
$$= O(\ln(n)^2 L(n)^{2\kappa} n^{-2\kappa(1-H)})$$
$$= o(1).$$

We choose $\rho$ large enough so that $(\rho - l)(1 - H) > 3$. Now, decompose $X_i$ as $X_i = \gamma_i U_i + \delta_i W_i$, where $EU_i^2 = EW_i^2 = 1$ and $\gamma_i U_i = E(X_i | \mathcal{F}_0)$. So,

$$E\left(V \sum_{k=\rho}^\infty g_k h_k(X_i) \Big| \mathcal{F}_0\right) = E\left(V \sum_{k=\rho}^\infty g_k \sum_{j=0}^k \binom{k}{j} \gamma_i^{k-j} \delta_i^j h_{k-j}(U_i) h_j(W_i) \Big| \mathcal{F}_0\right)$$
$$= \sum_{k=\rho}^\infty g_k \sum_{j=0}^l \binom{k}{j} \gamma_i^{k-j} \delta_i^j h_{k-j}(U_i) E(V h_j(W_i) | \mathcal{F}_0).$$

By [8], Theorem 5.10, $E(V h_j(W_i))^2$ is uniformly bounded over $i$ and $j$. Since $\gamma_i^2 = \sum_{j=i}^\infty b_j^2 \sim (2-2H) i^{2H-2} L_1(i)^2$, we have

$$\sup_{i > n^{1/3}, k \geq \rho, 0 \leq j \leq l} \binom{k}{j} \gamma_i^{k-j} = o(n^{-1})$$

and so when $i > n^{1/3}$,

$$n^2 E\left(\sum_{k=\rho}^\infty g_k \sum_{j=0}^l \binom{k}{j} \gamma_i^{k-j} \delta_i^j h_{k-j}(U_i)\right)^2 \to 0.$$

And, since $EZ_i^2 = O(n^{-1})$, it follows that

$$E\left(E\left(V\sum_{k=\rho}^{\infty} g_k h_k(X_i)\Big|\mathcal{F}_0\right)\right)^2 \to 0$$

as $n \to \infty$, which proves equation (2.2). By [8], Corollary 5.7,

$$E[Z_i Z_j] \le [EZ_i]^{2/(1+|EX_i X_j|)}$$

and when $|i-j|$ is large, $EX_i X_j > 0$, so $E[Z_i Z_j] = \sum_{i=0}^{\infty} g_k^2 k!(EX_i X_j)^k \ge g_0^2 = [EZ_i]^2$, from which equation (2.3) follows. □

PROPOSITION 2.1. *Let $\Delta^*$ be a finite union of finite rectangles in $(\mathbb{R}\setminus(-c, c)) \times \mathbb{R}$. Let $V$ be a random variable, measurable in the $\sigma$-algebra generated by $B(t)$ with $|V| \le 1$. Then,*

$$E[VI(\nu_n(\Delta^*) = 0)] \to EVE[I(\nu(\Delta^*) = 0)].$$

PROOF. For $\epsilon > 0$, we can partition $\Delta^* = \bigcup_{j=1}^{k} \Delta_j^* = \bigcup_{j=1}^{k} \Delta_j \times [s_j, t_j)$, where $\Delta_j$ are finite unions of finite intervals in $\mathbb{R}\setminus(-c,c)$, $s_1 < t_1 \le s_2 < \cdots \le t_k$ and $\sum_{j=1}^{k} \frac{1}{2}\mu_{\Delta_j}^2 [t_j - s_j]^2 < \epsilon$. We can write

$$I(\nu_n(\Delta_j^*) = 0) = \left(1 - E\sum_{i=\lfloor ns_j\rfloor+1}^{\lfloor nt_j\rfloor} Z_i\right) - \left(\sum_{i=\lfloor ns_j\rfloor+1}^{\lfloor nt_j\rfloor} (Z_i - EZ_i)\right)$$

$$+ \left(I(\nu_n(\Delta_j^*) = 0) - 1 + \sum_{i=\lfloor ns_j\rfloor+1}^{\lfloor nt_j\rfloor} Z_i\right).$$

By equation (2.4),

$$(2.5) \qquad \lim_n \prod_{j=1}^{k}\left(1 - E\sum_{i=\lfloor ns_j\rfloor+1}^{\lfloor nt_j\rfloor} Z_i\right) = \prod_{j=1}^{k}(1 - \mu_{\Delta_l}[t_l - s_l]).$$

By Lemma 2.1,

$$(2.6) \qquad E\left|V\left(\prod_{j=1}^{l-1} I(\nu_n(\Delta_l^*) = 0)\right)\left(\sum_{i=\lfloor ns_l\rfloor+1}^{\lfloor nt_l\rfloor} Z_i - EZ_i\right)\right| \to 0$$

as $n \to \infty$. Now, $\sum_{i=\lfloor ns_j\rfloor+1}^{\lfloor nt_j\rfloor} \sum_{j=i+1}^{\lfloor nt_j\rfloor} Z_i Z_j$ is the number of pairs of points in $\nu_n(\Delta_j^*)$ and so is equal to $\frac{1}{2}\nu_n(\Delta_j^*)(\nu_n(\Delta_j^*) - 1)$. This is greater than or equal to $\nu_n(\Delta_j^*) - 1$ when $\nu_n(\Delta_j^*) \ge 1$, so

$$0 \le I(\nu_n(\Delta_j^*) = 0) - 1 + \sum_{i=\lfloor ns_j\rfloor+1}^{\lfloor nt_j\rfloor} Z_i \le \sum_{i=\lfloor ns_j\rfloor+1}^{\lfloor nt_j\rfloor} \sum_{j=i+1}^{\lfloor nt_j\rfloor} Z_i Z_j$$

and so, by Lemma 2.1,

$$
\limsup_n E\left|V\left(\prod_{j=1}^{l-1} I(\nu_n(\Delta_l^*)=0)\right)\left(I(\nu_n(\Delta_l^*)=0)-1+E\sum_{i=\lfloor ns_l\rfloor+1}^{\lfloor nt_l\rfloor} Z_i\right)\right| \quad (2.7)
$$
$$
\leq \tfrac{1}{2}\mu_{\Delta_l}^2[t_l-s_l]^2.
$$

Putting together equations (2.5), (2.6) and (2.7), we obtain

$$
\limsup_n \left| E[VI(\nu_n(\Delta^*)=0)] - EV\prod_{j=1}^k (1-\mu_{\Delta_l}[t_l-s_l]) \right| \leq \epsilon.
$$

Taking increasingly fine partitions of $\Delta^*$, $\prod_{j=1}^k(1-\mu_{\Delta_l}[t_l-s_l]) = \prod_{j=1}^k(1-\nu'(\Delta_j^*))$ converges to $\exp(-\nu'(\Delta_j^*)) = EI(\nu(\Delta^*)=0)$, which completes the proof. □

Now, if $\Delta^*$ is a finite union of finite rectangles in $(\mathbb{R}\setminus(-c,c))\times\mathbb{R}$, then it can be written as $\Delta^* = \bigcup_{i=1}^k \Delta_i \times [s_i,t_i)$, where $\Delta_i$ are finite unions of finite intervals in $\mathbb{R}\setminus(-c,c)$ and $s_1 < t_1 \leq s_2 < \cdots \leq t_k$. Then, by Lemma 2.1, $E\nu_n(\Delta^*) \to \sum_{i=1}^k E\nu(\Delta^*)$ and, by Proposition 2.1, $EI(\nu_n(\Delta^*)=0) \to EI(\nu(\Delta^*)=0)$ as $n\to\infty$. Also, $\lim_{x\to\infty} \limsup_n P(\nu_n(\mathbb{R}\setminus(-x,x)\times[0,t])) = 0$, so, by [9] Theorems 4.7 and 4.9, $\nu_n$ converges weakly to $\nu$.

Now, suppose that $V$ is the indicator function of some event generated by $R_{\kappa,H}$. Then, by Lemma 2.1, $EV\nu_n(\Delta^*) \to \sum_{i=1}^k EVE\nu(\Delta^*)$ and, by Proposition 2.1, $P(\nu_n(\Delta^*)=0) \to P(\nu(\Delta^*)=0)$ as $n\to\infty$. So, asymptotically, this convergence takes place independently of $V$ and we can conclude that

$$
(2.8) \quad \left(L(n)^{-\kappa}n^{-(1-\kappa(1-H))}\sum_{i=1}^{\lfloor nt\rfloor} E(f(X_i)|\mathcal{F}_{i-d}),\nu_n\right) \to (f_\kappa R_{\kappa,H}(t),\nu)
$$

weakly, jointly where $R_{\kappa,H}$ and $\nu$ are independent.

By Karamata's lemma (see [10]),

$$
E(f(X_i)I(|f(X_i)|<ca_n))^2 \sim \frac{\alpha-1}{2-\alpha}c^{2-\alpha}a_n^2 n^{-1}
$$

as $n\to\infty$. If $Y_i = f(X_i)I(|f(X_i)|<ca_n) - E(f(X_i)I(|f(X_i)|<ca_n)|\mathcal{F}_{i-d})$, then $Y_k, Y_{k+d}, Y_{k+2d},\ldots$ are martingale differences for each $0\leq k<d$. Then, using Doob's maximal inequality,

$$
(2.9) \quad \limsup_n E\left(\max_{0\leq t\leq T} a_n^{-1}\sum_{i=1}^{\lfloor nt\rfloor} Y_i\right)^2 \leq 4dt^2\frac{\alpha-1}{2-\alpha}c^{2-\alpha}.
$$

Now, define $Y_i' = f(X_i)I(|f(X_i)|\geq ca_n)$. Again, by Karamata's lemma,

$$
E|Y_i'|^{1+\theta} \sim \frac{\alpha}{1+\theta}(ca_n)^{1+\theta-\alpha}L_3(ca_n)
$$

as $n \to \infty$. Applying Theorem 5.1 of [8], we have

$$E\left(a_n^{-1} \sum_{i=1}^{\lfloor nt \rfloor} E(Y_i'|\mathcal{F}_{i-d}) - EY_i'\right)^2$$
$$\leq C(nt)^{2H} L_1(nt)^2 a_n^{-2\alpha/(1+\theta)} L_3(ca_n)^{2/(1+\theta)} \to 0,$$

which implies tightness and finite-dimensional convergence, so

$$(2.10) \qquad a_n^{-1} \sum_{i=1}^{\lfloor nt \rfloor} (E(Y_i'|\mathcal{F}_{i-d}) - EY_i') \to 0$$

weakly.

Now, noting that

$$\lim_{x \to \infty} \limsup_n P(\nu_n(\mathbb{R} \setminus (-x, x) \times [0, t])) = 0,$$

it follows that

$$(2.11) \qquad a_n^{-1} \sum_{i=1}^{\lfloor nt \rfloor} f(X_i) I(a_n^{-1}|f(X_i)| \geq c) \to \int_{\mathbb{R} \setminus (-c,c)} \int_0^t y \, d\nu$$

weakly in the Skorohod topology. By the representation given in [11], Section 3.12,

$$(2.12) \qquad \int_{\mathbb{R} \setminus (-c,c)} \int_0^t y \, d(\nu - \nu') \to R^*(t)$$

in probability in the Skorohod $J_1$-topology as $c \to \infty$. Applying Theorem 4.2 of [1] to equations (2.9), (2.10), (2.11) and (2.12), we conclude that

$$(2.13) \qquad a_n^{-1} \sum_{i=1}^{\lfloor nt \rfloor} f(X_i) - E(f(X_i)|\mathcal{F}_{i-d}) \to R^*(t)$$

weakly as $n \to \infty$. Because of the asymptotic independence of $R_{\kappa, 1-\kappa(1-H)}$ and $\nu$,

$$\left(L(n)^{-\kappa} n^{-(1-\kappa(1-H))} \sum_{i=1}^{\lfloor nt \rfloor} E(f(X_i)|\mathcal{F}_{i-d}), a_n^{-1} \sum_{i=1}^{\lfloor nt \rfloor} f(X_i) - E(f(X_i)|\mathcal{F}_{i-d})\right)$$
$$\xrightarrow{d} (f_\kappa R_{\kappa, H}(t), R^*(t)),$$

where $R_{\kappa, 1-\kappa(1-H)}(t)$ and $R^*(t)$ are independent. Theorem 1.1 follows immediately.

The proof of Theorem 1.2 proceeds similarly. The measure $\nu_n$ can be constructed as before. To apply Theorem 4.2 of [1] and complete the proof, we need only show that

$$(2.14) \qquad \lim_{c \to 0} \limsup_n E\left(\max_{0 \le t \le T} a_n^{-1} \sum_{i=1}^{\lfloor nt \rfloor} Y_i - EY_i\right)^2 = 0,$$

where $Y_i = f(X_i)I(|f(X_i)| < ca_n)$, which follows from showing that the lower order chaos terms are insignificant in the limit.

DEPARTMENT OF STATISTICS
UNIVERSITY OF CALIFORNIA, BERKELEY
BERKELEY, CALIFORNIA 94720
USA
E-MAIL: sly@stat.berkeley.edu

CENTRE FOR MATHEMATICS AND ITS APPLICATIONS
MATHEMATICAL SCIENCES INSTITUTE
AUSTRALIAN NATIONAL UNIVERSITY
ACT 0200
AUSTRALIA
E-MAIL: Chris.Heyde@maths.anu.edu.au